\documentclass[10pt,twoside]{article}
\usepackage{graphicx}
\usepackage{amsmath,amsthm,amsfonts,verbatim}
\usepackage{Latex-document}

\theoremstyle{plain}
 \newtheorem{theorem}{Theorem}

 \newtheorem{problem}{Problem}
 \newtheorem{open}{Open Problem}
 
 \theoremstyle{definition}

\newcommand{\distortion}{\mbox{distortion}}
\newcommand{\contraction}{\mbox{contraction}}
\newcommand{\expansion}{\mbox{expansion}}

\title{\bf Finite Metric Spaces---Combinatorics,\vskip -2mm Geometry and Algorithms\vskip 6mm}

\author{Nathan Linial\thanks{School of Computer Science and Engineering, Hebrew University, Jerusalem
91904, Israel. E-mail:nati@cs.huji.ac.il}\vspace*{-0.5cm}}

\date{\vspace{-8mm}}

\begin{document}

\maketitle

\thispagestyle{first} \setcounter{page}{573}

\begin{abstract}\vskip 3mm

Finite metric spaces arise in many different contexts. Enormous bodies of data, scientific, commercial and others
can often be viewed as large metric spaces. It turns out that the metric of graphs reveals a lot of interesting
information. Metric spaces also come up in many recent advances in the theory of algorithms. Finally, finite
submetrics of classical geometric objects such as normed spaces or manifolds reflect many important properties of
the underlying structure. In this paper we review some of the recent advances in this area.

\vskip 4.5mm

\noindent {\bf 2000 Mathematics Subject Classification:} Combinatorics, Algorithms, Geometry.

\end{abstract}

\vskip 12mm

\section{Introduction}

\vskip-5mm \hspace{5mm}

The constantly intensifying ties between combinatorics and
geometry are among the most significant developments in Discrete
Mathematics in recent years. These connections are manifold, and
it is, perhaps, still too early to fully evaluate this
relationship. This article deals only with what might be called
{\em the geometrization of combinatorics}. Namely, the idea that
viewing combinatorial objects from a geometric perspective often
yields unexpected insights. Even more concretely, we concentrate
on finite metric spaces and their embeddings.

To illustrate the underlying idea, it may be best to begin with a
practical problem. There are many disciplines, scientific,
technological, economic and others, which crucially depend on the
analysis of large bodies of data. Technological advances have made
it possible to collect enormous amounts of interesting data, and
further progress depends on our ability to organize and classify
these data so as to allow meaningful and insightful analysis. A
case in point is bioinformatics where huge bodies of data - DNA
sequences, protein sequences, information about expression levels
etc. all await analysis. Let us consider, for example, the space
of all proteins. For the purpose of the current discussion, a
protein may be viewed as a word in an alphabet of 20 letters
(amino acids). Word lengths vary from under fifty to several
thousands, the most typical length being several hundred letters.
At this writing, there are about half a million proteins whose
sequence is known. Algorithms were developed over the years to
evaluate the similarity of different proteins, and there are
standard computer programs that calculate distances among proteins
very efficiently. This turns the collection of all known proteins
into a metric space of about half a million elements. Proper
analysis of this space is of great importance for the biological
sciences. Thus, this huge body of sequence data takes a geometric
form, namely, a finite metric space, and it becomes feasible to
use geometric concepts and tools in the analysis of this data.

In the combinatorial realm proper, and in the design and analysis
of algorithms, similar ideas have proved very useful as well. A
graph is completely characterized by its (shortest path, or
geodesic) metric. The analysis of this metric provides a lot of
useful information about the graph. Moreover, given a graph $G$,
one may modify $G$'s metric by assigning nonnegative {\em lengths}
to $G$'s edges. By varying these edge lengths, a family of finite
metrics is obtained, the properties of which reflect a good deal
of structural information about $G$. We mention in passing that
there are other useful and interesting geometric viewpoints of
graphs. Thus, it is useful to geometrically realize a graph by
assigning vectors to the vertices and posit that adjacent vertices
correspond to orthogonal vectors. Graphs can encode the
intersection patterns of geometric objects. These are all
interesting instances of our basic paradigm: In the study of
combinatorial objects, and especially graphs, it is often
beneficial to develop a perspective from which the graph is
perceived geometrically.

Aside from what has already been thus accomplished, this approach
holds a great promise. Combinatorics as we know it, is still a
very young subject. (There is no official date of birth, and Euler
was undoubtedly a giant in our field, but I think that the dawn of
modern combinatorics can be dated to the 1930's). Discrete
Mathematics stands to gain a lot from interactions with older,
better established fields. This geometrization of combinatorics
indeed creates clear and tangible connections with various
subfields of geometry. So far the study of finite metric spaces
has had substantial connections with the theory of
finite-dimensional normed spaces, but it seems safe to predict
that useful ties with differential geometry will soon emerge. With
the possible incorporation of probabilistic tools, now commonplace
in combinatorics, we can expect very exciting outcomes.

A good sign for the vitality of this area is the large number of
intriguing open problems. We will present here some of those that
we particularly like. In a recent meeting (Haifa, March '02), a
list of open problems in this area has been collected, see
http://www.kam.mff.cuni.cz/$\tilde{}$matousek/haifaop.ps. More
extensive surveys of this area can be found
in~\cite{matousek_disc_geom} Chapter 15, and~\cite{Ind_surv}.

In view of this description, it should not come as a surprise to
the reader that this theory is characterized as being
\begin{itemize}
\item
{\em Asymptotic:} We are mostly interested in analyzing large,
finite metric spaces, graphs and data sets.
\item
{\em Approximate:} While it is possible to postulate that the
geometric situation agrees perfectly with the combinatorics, it is
much more beneficial to investigate the approximate version. This
leads to a richer theory that is quantitative in nature. Rather
than a binary question whether perfect mimicking is possible or
not, we ask {\em how well} a given combinatorial object can be
approximated geometrically.
\item
{\em Algorithmic:} Existential results are very important and
interesting in this area, but we always prefer it when such a
result is accompanied by an efficient algorithm.
\item
It is mostly {\em comparative:} There are certain classes of
finite metric spaces that we favor. These may have a particularly
simple structure or be very well understood. Other, less well
behaved spaces are being compared to, and approximated by, these
``nice" metrics.
\end{itemize}

So, how should we compare between two metrics? Let $(X,d)$ and
$(Y,\rho)$ be two metric spaces and let $\varphi: X \rightarrow Y$
be a mapping between them. We quantify the extent to which
$\varphi$ expands, resp. contracts distances: $\expansion(\varphi)
= \sup_{x,y \in X} \frac{\rho(\varphi(x),\varphi(y))}{d(x,y)}$ and
$\contraction(\varphi) = \sup_{x,y \in X}
\frac{d(x,y)}{\rho(\varphi(x),\varphi(y))}$.

Finally, the main definition is: $\distortion(\varphi) =
\expansion(\varphi) \cdot \contraction(\varphi)$.

In other words, we consider the tightest constants $\alpha \ge
\beta$ for which $\alpha \ge
\frac{\rho(\varphi(x),\varphi(y))}{d(x,y)} \ge \beta$ always
holds, and define $\distortion(\varphi)$ as
$\frac{\alpha}{\beta}$. We call $\varphi$ an {\em isometry} when
$\distortion(\varphi) = 1$. This deviates somewhat from the
conventional definition, and a map that multiplies all distances
by a constant (not necessarily $1$) is being considered here as an
isometry.

The least distortion with which $(X,d)$ can be embedded in
$(Y,\rho)$ is denoted $c_Y(X)=c_Y(X,d)$. If $\cal C$ is a class of
metric spaces, then the infimum of $c_Y(X)$ over all $Y \in \cal
C$ is denoted by $c_{\cal C}(X)$. When $\cal C$ is the class of
finite-dimensional $l_p$ spaces $\{l_p^n | n=1,2,\dots\}$ we
denote $c_{\cal C}(X)$ by $c_p(X)$.

One of the {\bf major problems} in this area is:
\begin{problem}
\rm Given a finite metric space $(X,d)$ and a class of metrics
$\cal C$, find the (nearly) best approximation for $X$ by a metric
from $\cal C$. In other words, find a metric space $Y \in \cal C$
and a map $\varphi: X \rightarrow Y$ such that
$\distortion(\varphi)$ (nearly) equals $c_{\cal C}(X)$.
\end{problem}

The classes of metric spaces $\cal C$ for which this problem has
so far been studied are: (i) Metrics of normed spaces, especially
$l_p^n$ for $\infty \ge p \ge 1$ and $n=1,2,\ldots$. (ii) Metrics
of special families of graphs, most notably trees, as well as
convex combinations thereof.

One more convention: Speaking of $l_p$, either means infinite
dimensional $l_p$, or, what is often the same, that we do not care
about the dimension of the space in which we embed a given metric.

To get a first feeling for this subject, let us consider the
smallest nontrivial example. Every $3$-point metric embeds
isometrically into the plane, but as we show now, the metric of
$K_{1,3}$, the $4$-vertex tree with a root and three leaves, has
no isometric embedding into $l_2$. Let $x$, resp. $y_i$ be the
image of the root and the leaves of this tree. Since $d(x,y_i)=1$
and $d(y_i,y_j)=2$ for all $i \ne j$, it follows that the three
points $x,y_i,y_j$ are colinear for every $i \ne j$. Thus, all
four points are colinear, leading to a contradiction. It can be
shown that the least distorted image of this graph in $l_2$ is in
the plane with $120^{\circ}$ degree angle among the edges. Below
(Section~\ref{sec_l2}) we present a polynomial-time algorithm that
determines $c_2(X)$, the least $l_2$ distortion for any finite
metric $(X,d)$.

Another easy fact which belongs into this warm-up section is that
$c_{\infty}(X) = 1$ for every finite metric $(X,d)$. That is, the
space $l_{\infty}$ space contains an isometric copy of every
finite
metric space.\\
{\bf Acknowledgment:} Helpful remarks on this article by R.
Krauthgamer, A. Magen, J. Matou{\v{s}}ek, and Yu. Rabinovich are
gratefully acknowledged.

\section{Embedding into {\boldmath $l_2$}} \label{sec_l2}

\vskip-5mm \hspace{5mm}

This is by far the most developed part of the
theory. There are several good reasons for this part of the theory
to have attracted the most attention so far. Consider the
practical context, where a metric space represents some large data
set, and where the major driving force is the search for good
algorithms for data analysis. If the data set you need to analyze
happens to be a large set of points in $l_2$, there are many tools
at your disposal, from geometry, algebra and analysis. So if your
data can be well approximated in $l_2$, this is of great practical
advantage. There is another reason for the special status of $l_2$
in this area. To explain it, we need to introduce some terminology
from Banach space theory. The {\em Banach-Mazur distance} among
two normed spaces $X$ and $Y$, is said to be $\le c$, if there is
a {\em linear} map $\varphi: X \rightarrow Y$ with
$\distortion(\varphi) \le c$. What we are doing here may very well
be described as a search for the metric counterpart of this highly
developed linear theory. See~\cite{MS86} for an introduction to
this field and~\cite{BeLin} for a comprehensive cover of the
nonlinear theory. The grandfather of the linear theory is the
celebrated theorem of Dvoretzky~\cite{dvor}.

\begin{theorem}[Dvoretzky]
\label{dvoretzky} For every $n$ and $\epsilon >0$, every
$n$-dimensional normed space contains a $k=\Omega(\epsilon^2 \cdot
\log n)$-dimensional space whose Banach-Mazur distance from $l_2$
is $\le 1+\epsilon$.
\end{theorem}

Thus, among embeddings into normed spaces, embeddings into $l_2$
are the hardest to come by.

We begin our story with an important theorem of
Bourgain~\cite{Bourgain85}.
\begin{theorem}
\label{bour_thm} Every $n$-point metric space \footnote{Here and
elsewhere, unless otherwise stated, $n = |X|$, the cardinality of
the metric space in question.} embeds in $l_2$ with distortion
$\le O(\log n)$.
\end{theorem}
Not only is this a fundamental result, Bourgain's proof of the
theorem readily translates into an efficient randomized algorithm
that finds, for any given finite $(X,d)$ an embedding in $l_2$ of
distortion $\le O(\log n)$. The algorithm is so simple that we
record it here. Given the metric space $(X,d)$, we map every point
$x \in X$ to $\varphi(x)$, an $O(\log^2 n)$-dimensional vector.
Coordinates in $\varphi(\cdot)$ correspond to subsets $S \subseteq
X$, and the $S$-th coordinate in $\varphi(x)$ is simply $d(x,S)$,
the minimum of $d(x,y)$ over all $y \in S$. To define the map
$\varphi$, we need to specify, then, the collection of subsets $S$
that we utilize. These sets are selected randomly. Namely, you
randomly select $O(\log n)$ sets of size $1$, another $O(\log n)$
sets of size $2$, of size $4, 8...,\frac{n}{2}$.

In view of Bourgain's Theorem, several questions suggest
themselves naturally:

\begin{itemize}
\item
Is this bound tight? The answer is positive, see
Theorem~\ref{expander}.
\item
Given that $\max c_2(X)$ over all $n$-point metrics is
$\Theta(\log n)$, what about metrics that are closer to $l_2$? Is
there a polynomial-time algorithm to compute $c_2(X,d)$ (That is,
the least distortion in an embedding of $X$ into $l_2$)? Again the
answer is affirmative, see below and Theorem~\ref{sdp}.
\item
Are there interesting families of metric spaces for which $c_2$ is
substantially smaller than $\log n$? Indeed, there are, see, e.g.,
Theorem~\ref{tree}.
\end{itemize}

So let us proceed with the answers to these questions. {\em
Expanders} are graphs which cannot be disconnected into two large
subgraphs by removing relatively few edges. Specifically, a graph
$G$ on $n$ vertices is said to be an $\epsilon$-{\em
(edge)-expander} if, for every set $S$ of $\le n/2$ vertices,
there are at least $\epsilon |S|$ edges between $S$ and its
complement. It is said to be $k$-{\em regular} if every vertex has
exactly $k$ neighbors. The theory of expander graphs is a
fascinating chapter in discrete mathematics and theoretical
computer science. It is not obvious that arbitrarily large
$k$-regular graphs exist with expansion $\epsilon$ bounded away
from zero. In fact, in the early days of this area, conjectures to
the contrary had been made. It turns out, however, that expanders
are rather ubiquitous. For every $k \ge 3$, the probability that a
randomly chosen $k$-regular graph has expansion $\epsilon > k/10$
tends to $1$ as the number of vertices $n$ tends to $\infty$. It
turns out that the metrics of expander graphs are as far from
$l_2$ as possible. \footnote{We freely interchange between a graph
and its (shortest path) metric.}

\begin{theorem}[\cite{LLR95}, see also~\cite{Matousek97,LM00}]
\label{expander} Let $G$ be an $n$-vertex $k$-regular
$\epsilon$-expander graph ($k \ge 3$, $\epsilon >0$). Then $c_2(G)
\ge c \log n$ where $c$ depends only on $k$ and $\epsilon$.
\end{theorem}

Metric geometry is by no means a new subject, and indeed metrics
that embed isometrically into $l_2$ were characterized long ago
(see e.g. \cite{Blumenthal}). This is a special case of the more
recent results. Let $\varphi: X \rightarrow l_2^n$ be an
embedding. The condition that $\distortion(\varphi) \le c$ can be
expressed as a system of linear inequalities in the entries of the
Gram matrix corresponding to the vectors in $\varphi(X)$.
Therefore, the computation of $c_2(X)$ is an instance of {\em
semidefinite quadratic programming} and can be found in polynomial
time. \footnote{This is not quite accurate. Given an $n$-point
space $(X,d)$ and $\epsilon > 0$, the algorithm can determine
$c_2(X,d)$ with relative error $<\epsilon$ in time polynomial in
$n$ and $\frac{1}{\epsilon}$.} This formulation of the problem
has, however, other useful consequences. The duality principle of
convex programming yields a max-min formula for $c_2$.

\begin{theorem}[\cite{LLR95}]
\label{sdp} For every finite metric space $(X,d)$,
$$
c_2(X,d) = \max \sqrt{{\frac{\sum_{i,j:q_{i,j} >
0}d^2(i,j)q_{i,j}}
  {\sum_{i,j:q_{i,j}  <  0}d^2(i,j)|q_{i,j}|}}},
$$
where the maximum is over all matrices $Q$ so that
\begin{enumerate}
\item
$Q$ is positive semidefinite, and
\item
The entries in every row in $Q$ sum to zero.
\end{enumerate}
\end{theorem}

Consider the metric of the $r$-dimensional cube. As shown by
Enflo~\cite{Enflo69}, the least distorted embedding of this metric
is simply the identity map into $l_2^r$, which has distortion
$\sqrt{r}$. Our first illustration for the power of the quadratic
programming method is that we provide a quick elementary proof for
this fact, earlier proofs of which required heavier machinery. The
rows and columns of the matrix $Q$ are indexed by the $2^r$
vertices of the $r$-dimensional cube. The $(x,y)$ entry of $Q$ is:
(i) $r-1$ if $x=y$, (ii) It is $-1$ if $x$ and $y$ are neighbors
(they are represented by two $0,1$ vectors that differ in exactly
one coordinate, and (iii) It is $1$ if $x$ and $y$ are antipodal,
i.e., they differ in all $r$ coordinates. (iv) All other entries
of $Q$ are zero. We leave out the details and only indicate how to
prove that $Q$ is positive semidefinite. It is possible to express
$Q=(r-1)I - A + P$, where $A$ is the adjacency matrix of the
$r$-cube and $P$ is the (permutation) matrix corresponding to
being antipodal. The eigenfunctions of $A$ are well known, namely,
they are the $2^r$ Walsh functions. The same vectors happen to be
also the eigenvectors of $Q$ and all have nonnegative eigenvalues.

As another application of this method (also from~\cite{LM00}),
here is a quick proof of Theorem~\ref{expander}. It is
known~\cite{Alon_eigen_expan} that if $G$ is a $k$-regular
$\epsilon$-expander graph and $A$ is $G$'s adjacency matrix, then
the second eigenvalue of $A$ is $<k - \delta$ for some $\delta$
that depends on $k$ and $\epsilon$, but not on the size of the
graph \footnote{$A$'s first eigenvalue is clearly $k$. This is the
combinatorial analogue of Cheeger's Theorem~\cite{cheeger} about
the spectrum of the Laplacian.}. It is not hard to show that the
vertices of a graph with bounded degrees can be paired up so that
every two paired vertices are at distance $\Omega(\log n)$. Let
$P$ be the permutation matrix corresponding to such a pairing. It
is not hard to establish Theorem~\ref{expander} using the matrix
$Q = kI - A + \frac{\delta}{2} (P-I)$. More sophisticated
applications of this method will be described below
(Theorem~\ref{girth_l2}).

\section{Specific families of graph metrics}

\vskip-5mm \hspace{5mm}

For various graph families, it is possible find embeddings into $l_2$ with distortion asymptotically smaller than
$\log n$. This often applies as well to graphs with arbitrary nonnegative edge lengths.

\subsection{Trees}
\vskip -5mm \hspace{5mm}

The metrics of trees are quite restricted. They can be
characterized through a four-term inequality
(e.g.~\cite{deza_lau}). It is also not hard to see that every tree
metric embeds isometrically into $l_1$. They can also be embedded
into $l_2$ with a relatively low distortion.
\begin{theorem}[Matou{\v{s}}ek~\cite{Matousek99}]
\label{tree} Every tree on $n$ vertices can be embedded into $l_2$
with distortion $\le O(\sqrt{\log \log n})$.
\end{theorem}
Bourgain~\cite{Bourgain86} had earlier shown that this bound is
attained for complete binary trees. (See~\cite{Li:Sa_tree} for an
elementary proof of this.)

\subsection{Planar graphs}
\vskip -5mm \hspace{5mm}

It turns out that the metrics of planar graphs have good embedding
into $l_2$. Rao~\cite{Rao99} showed:
\begin{theorem}
\label{rao_thm} Every planar graph embeds in $l_2$ with distortion
$O(\sqrt{\log n})$.
\end{theorem}
A recent construction of Newman and
Rabinovich~\cite{RabinovichNewmanSoCG} shows that this bound is
tight.

\subsection{Graphs of high girth}\vskip -5mm \hspace{5mm}

The {\em girth} of a graph is the length of the shortest cycle in
the graph. If you restrict your attention (as we do in this
section) to graphs in which all vertex degrees are $\ge 3$, then
it is still a major challenge to construct graphs with very high
girth, i.e., having no short cycles. The metrics of such graphs
seem far from $l_2$, so in~\cite{LLR95} it was conjectured that
$c_2(G) \ge \Omega(g)$ for every graph $G$ of girth $g$ in which
all vertex degrees are $\ge 3$. There are known examples of
$n$-vertex $k$-regular expanders whose girth is $\Omega(\log n)$.
In view of Theorem~\ref{bour_thm}, such graphs show that this
conjecture, if true, is best possible. Recently, the following was
shown:
\begin{theorem}[\cite{LiMaNa_girth}]
\label{girth_l2} Let $G$ be a $k$-regular graph $k \ge 3$ with
girth $g$. Then $c_2(G) \ge \Omega(\sqrt{g})$.
\end{theorem}

Two proofs of this theorem are given in~\cite{LiMaNa_girth}. One
is based on the notion of {\em Markov Type} due to
Ball~\cite{ball_markov_type}. The underlying idea of this proof is
that a random walk on a graph with girth $g$ and all vertex degree
$\ge 3$ drifts at a constant speed away from its starting point
for time $\Omega(g)$. On the other hand, in an appropriately
defined class of random walks in Euclidean space, at time $T$ the
walk is expected to be only $O(\sqrt{T})$ away from its origin. If
we compare between the graph itself and its image under an
embedding in $l_2$, this discrepancy must be accounted for by a
metrical distortion. The comparison at time $T = \Theta(g)$ yields
a distortion of $\Omega(\sqrt{g})$.

The other proof again employs semidefinite programming, using the
matrix $Q=\alpha I - A + \beta B$. Here $A$ is the graph's
adjacency matrix, and $B$ is a $0,1$ matrix where $B_{xy} = 1$ if
$x$ and $y$ are at distance $g/2$ in $G$. The parameters $\alpha$
and $\beta$ have to satisfy the two conditions from
Theorem~\ref{sdp}. A key observation is that due to the high
girth, $B$ can be expressed as $P_{g/2}(A)$ where $P_j$ is the
$j$-th {\em Geronimus Polynomial}, a known family of orthogonal
polynomials. The proof depends on the distribution of zeros for
these polynomials, and other analytical properties that they have.

Our present state of knowledge leads us to ask:

\begin{open}
\rm How small can $c_2(G)$ be for a a graph $G$ of girth $g$ in
which all vertices have degree $\ge 3$? The answer lies between
$\Omega(\sqrt{g})$ and $O(g)$.
\end{open}

An earlier result of Rabinovich and Raz~\cite{Rab:Raz} reveals
another connection between high girth and distortion. Let
$\varphi$ be a map from a graph of girth $g$ to a graph of smaller
{\em Euler characteristic} ($|E|-|V|+1$). Then
$\distortion(\varphi) \ge \Omega(g)$.

\section{Algorithmic applications} \label{alg_sec}

\vskip-5mm \hspace{5mm}

Among the most pleasing aspects of this field, are
the many beautiful applications it has to the design of new
algorithms.

\subsection{Multicommodity flow and sparsest cuts}\vskip -5mm \hspace{5mm}

Flows in networks are a classical subject in discrete optimization
and a topic of many investigations (see~\cite{schrijver} for a
comprehensive coverage). You are given a {\em network} i.e., a
graph with two specified vertices: The {\em source} $s$ and the
{\em sink} $t$. Edges have nonnegative {\em capacities}. The
objective is to ship as much of a given commodity between $s$ and
$t$, subject to two conditions: (i) In every vertex other than $s$
and $t$, matter is conserved, (ii) The flow through any edge must
not exceed the edge capacity. Let the set $S$ {\em separate} the
vertices $s$ and $t$, i.e., it contains exactly one of them.
Define $S$'s {\em capacity} as the sum of edge capacities over
those edges that connect $S$ to its complement. The {\em Max-flow
Min-cut} Theorem states that the largest possible flow equals the
minimum such capacity.

Here we consider the $k$-{\em commodity} version: Now there are
$k$ source-sink pairs $s_i, t_i, i=1,2,...,k$ for the $i$-th
commodity, and the $i$-th {\em demand} is $D_i >0$. We seek to
determine the largest $\phi > 0$ for which it is possible to flow
$\phi \cdot D_i$ of the $i$-th commodity between $s_i$ and $t_i$,
simultaneously for all $k \ge i \ge 1$ subject to conditions (i)
and (ii) above where in (ii) the {\em total} flow through an edge
should not exceed its capacity. With every subset of the vertices
$S$ we associate $\gamma(S)=\frac{\mbox{cap}(S)}{\mbox{dem}(S)}$.
As before, $\mbox{cap}(S)$ is the sum of the capacities of edges
between $S$ and its complement. The denominator $\mbox{dem}(S)$ is
$\sum D_i$ over all indices  $i$ so that $S$ separates $s_i$ and
$t_i$. It is trivially true that $\phi \le \gamma(S)$, for every
flow and every set $S$, but unlike the one-commodity case, $\min
\gamma(S)$ (the {\em sparsest cut}) need not equal $\max \phi$. As
for the algorithmic perspective, finding $\max \phi$ is a linear
program, so it can be computed in polynomial time. However, it is
$NP$-hard to determine the sparsest cut. Also, it is interesting
to find out how far $\max \phi$ and $\min \gamma(S)$ can be.
Consider the case where the underlying graph is an expander, edges
have unit capacities and every pair of vertices form a source-sink
pair with a unit demand. It is not hard to see that in this case
$\phi \le O(\frac{\min \gamma(S)}{\log n})$. On the other hand,

\begin{theorem}[\cite{LLR95}, see also \cite{AuRa:SparseCut}]
\label{multicommodity} In the $k$-commodity problem
$$
\max \phi \ge \Omega(\frac{\min \gamma(S)}{\log k}).
$$
\end{theorem}
We will be able to review the proof in Section~\ref{section_l1}.

\subsection{Graph bandwidth}\vskip -5mm \hspace{5mm}

In this computational problem, we are presented with an $n$-vertex
graph $G$. It is required to label the vertices with distinct
labels from $\{1,\ldots,n\}$ so that the difference between the
labels of any two adjacent vertices is not too big. Namely,
$$
\mbox{bw}(G) = \min_{\psi} \max_{xy \in E(G)} |\psi(x) - \psi(y)|,
$$
where the minimum is over all $1:1$ maps $\psi: V \rightarrow
\{1,\ldots,n\}$.

It is $NP$-hard to compute this parameter, and for many years no
decent approximation algorithm was known. However, a recent paper
by Feige~\cite{Feige00} provides a polylogarithmic approximation
for the bandwidth. The statement of his algorithm is simple enough
to be recorded here:
\begin{enumerate}
\item
Compute (a slight modification of) the embedding $\varphi: G
\rightarrow l_2$ that appears in the proof of Bourgain's
Theorem~\ref{bour_thm}.
\item
Select a random line $l$ and project $\varphi(G)$ onto it.
\item
Label the vertices of $G$ by the order at which their images
appear along the line $l$.
\end{enumerate}
Let $\beta(G) := \max_{x,r} \frac{|B_r(x)|}{r}$ where $B_r(x)$ is
the set of those vertices in $G$ at distance $\le r$ from $x$.
It's easy to see that $\mbox{bw}(G) \ge \Omega(\beta(G))$ and an
interesting feature of Feige's proof is that it shows that
$\mbox{bw}(G) \le O(\beta(G)\log^c n)$. His paper gives $c=3.5$
which was later~\cite{vemp_dun} improved to $c=3$.
\begin{open}
\rm Is it true that $\mbox{bw}(G) \le O(\beta(G)\log n)$?
\end{open}
It is not hard to see that this bound would be tight for
expanders.

\subsection{Bartal's method} \label{sect_bartal}

\vskip-5mm \hspace{5mm}

The following general structure theorem of
Bartal~\cite{bartal_trees} has numerous algorithmic applications:

\begin{theorem}
\label{bartal} For every finite metric space $(X,d)$ there is a
collection of trees $\{T_i~|~i \in I\}$, each of which has $X$ as
its set of leaves, and positive weights $\{p_i~|~i \in I\}$ with
$\sum_I p_i = 1$. Each of these tree metrics {\em dominates} $d$,
i.e., $dist_{T_i}(x,y) \ge d(x,y)$ for every $i$ and every $x,y
\in X$. On the other hand, for every $x,y \in X$,
$$
\sum_i p_i \cdot dist_{T_i}(x,y) \le O(\log n \cdot \log \log n
\cdot d(x,y)).
$$
\end{theorem}

{\sl Bartal's algorithmic paradigm} is a general principle
underlying the numerous algorithmic applications of this theorem:
Given an algorithmic problem on input a graph or a general metric
space $(X,d)$, find a collection of tree metrics $T_i$ and weights
$p_i$ as in Theorem~\ref{bartal}. Select one of the trees at
random, where $T_i$ is selected with probability $p_i$. Now solve
the problem for input $T_i$. (This description assumes, and this
is often the case, that the original optimization problem is
$NP$-hard in general, but feasible for tree metrics.

There are two features of the proof that we'd like to mention:\\
The trees $T_i$ are {\em HST}'s. In such trees, edge lengths
decrease exponentially as you move from the root toward the
leaves. They feature prominently in many
recent developments in this area.\\
The proof makes substantial use of {\em sparse decompositions} of
graphs. Given a graph, one seeks a probability distribution on all
partitions of the vertex set, so that (i) Parts have small
diameters (ii) Adjacent vertices are very likely to reside in the
same part. Such partitions have proved instrumental in the design
of many algorithms. In fact, an important tool in Rao's
Theorem~\ref{rao_thm} was an earlier result~\cite{Kl:Pl:Rao} about
the existence of very sparse partitions for the members of any
minor-closed families family of graphs.

\section{The mysterious {\boldmath $l_1$}} \label{section_l1}

\vskip-5mm \hspace{5mm}

We know much less about metric embeddings into
$l_1$, and the attempts to understand them give rise to many
intriguing open problems. We start by defining the {\em cut
metric} $d_S$ on $X$ where $S \subseteq X$, as follows: $d_S(x,y)
= 1$ if $x,y$ are separated by $S$ and is zero otherwise. A
simple, but useful observation is that the collection of all
$n$-point metrics in $l_1$ form a cone $\cal C$ whose extreme rays
are the cut metrics. \footnote{For each $n$, the $n$-point metrics
in $l_1$ form a cone $\cal C \rm_n$, but we suppress the index
$n$.} The book~\cite{deza_lau} provides a coverage of this area.

We are now able to complete the proof of
Theorem~\ref{multicommodity}. We retain the terminology of the
discussion around that theorem. Linear programming duality yields
the following alternative expression for the maximum $k$-commodity
flow problem on $G=(V,E)$:
$$
\label{dual_multi} \max \phi = \min \frac{\sum_E d(i,j) \cdot
c_{ij}}{\sum_1^k D_j \cdot d(s_j,t_j)}.
$$
Here the minimum is over all graphical metrics $d$ on $G$. Namely,
you assign nonnegative {\em lengths} to $G$'s edges and $d$ is the
induced shortest path metric on $G$'s vertices. Now let $d$ be the
graphical metric that minimizes this expression. A slight
adaptation of Bourgain's embedding algorithm yields an $l_1$
metric $\rho$ so that $\rho(i,j) \le d(i,j)$ for all $i,j$ and
$\rho(s_j,t_j) \ge \Omega(\frac{d(s_j,t_j)}{\log k})$ for all $j$.
But the minimum of $\frac{\sum_E \rho(i,j) \cdot c_{ij}}{\sum_1^k
D_j \cdot \rho(s_j,t_j)}$ over $l_1$ metrics is attained for
$\rho$ a cut metric, since cut metrics are the extreme rays of the
cone of $l_1$ metrics $\cal C$. This minimum, over cut metrics is
simply $\min \gamma(S)$, the sparsest cut value of the network.
The conclusion follows.

The identification between $l_1$ metrics and the cut cone $\cal C$
makes it desirable to find an algorithm to solve linear
optimization problems whose feasible set is this convex cone. Such
an algorithm would solve at one fell swoop a host of interesting
(and hard) problems such as max-cut, graph bisection and more.
This hope is hard to realize, since the ellipsoid method
(e.g.~\cite{schrijver}) applies only to convex bodies for which we
have efficient {\em membership and separation oracles}. For the
convex cone $\cal C$, that would mean that we need to efficiently
determine whether a given a real symmetric matrix $M$, represents
the metric on $n$ points in $l_1$. Moreover, if not, we ought to
find a hyperplane (in $n^2$ dimensions) that separates $M$ from
$\cal C$. Unfortunately, these questions are $NP$-hard
(e.g.~\cite{deza_lau}). It becomes, therefore, interesting to {\em
approximate} the cone $\cal C$. So, can we find another cone that
is close to $\cal C$ and for which computationally efficient
membership and separation oracles exist? There is a natural
candidate for the job. We say that a matrix $M$ is in {\em
square}-$l_2$, if there are points $x_i$ in $l_2$ such that
$M_{ij} = \|x_i - x_j\|_2^2$. Let $\cal S$ be the collection of
all all square-$l_2$ matrices which are also a metric (i.e. the
entries in $M$ also satisfy the triangle inequality). It is not
hard to see that $\cal C \subseteq \cal S$, but we ask:
\begin{open}
\rm What is the smallest $\alpha = \alpha(n)$, such that every $n
\times n$ matrix $M \in \cal S$ can be embedded in $l_1$ with
distortion $\le \alpha$ ?
\end{open}

It is not hard to see that every finite $l_2$ metric embeds
isometrically into $l_1$. But what about the opposite direction?
\begin{open}
\label{l1_to_l2}\rm Find $\max c_2(X)$ over all $(X,d)$ that are
$n$-point metrics in $l_1$. As we saw above, for the $n=2^r$
vertices of the $r$-cube the answer is $\sqrt{r} = \sqrt{\log n}$.
We suspect that this is the extreme case. No example is known
where $c_2$ is asymptotically larger that $\sqrt{\log n}$.
\end{open}

\subsection{Dimension reduction}\vskip -5mm \hspace{5mm}

Let us return to the applied aspect of this area. Even when a
given metric space can be approximated well in some normed space,
the {\em dimension} of the host space is quite significant. Data
analysis and clustering in $l_2^N$ for large $N$ is by no means
easy. In fact, practitioners in these areas often speak about the
{\em curse of dimensionality} when they refer to this problem. In
$l_2$ there is a basic result that answers this problem.
\begin{theorem}[Johnson Lindenstrauss~\cite{JL84}]
Every $n$-point metric in $l_2$ can be embedded into $l_2^k$ with
distortion $< 1+\epsilon$ where $k \le O(\frac{\log
n}{\epsilon^2})$.
\end{theorem}
Here, again, the proof yields an efficient randomized algorithm.
Namely, select a random $k$-dimensional subspace and project the
points to it.

What is the appropriate analogue of this theorem for $l_1$
metrics?
\begin{open}\rm
What is the smallest $k = k(n,\epsilon)$ so that every $n$-point
metric in $l_1$ can be embedded into $l_1^k$ with distortion $<
1+\epsilon$?
\end{open}
We know very little at the moment, namely $\Omega(\log n) \le k
\le O(n \log n)$ for constant $\epsilon >0$. The lower bound is
trivial and the upper bound is from~\cite{Sch_l1,talagrand_l1}.
Note that if the truth is at the lower bound, then this provides
an affirmative answer to Open Problem~\ref{l1_to_l2}.

\subsection{Planar graphs and other minor-closed families}
\vskip -5mm \hspace{5mm}

One of the most fascinating problems about $l_1$ metrics is:
\begin{open}\rm
Is there is an absolute constant $C > 0$ so that every metric of a
planar graph embeds into $l_1$ with distortion $< C$?
\end{open}
Even more daringly, the same can be asked for every minor-closed
family of graphs. Some initial success for smaller graph families
has been achieved already~\cite{Gup:New:Rab:Sin}.

\subsection{Large girth}\vskip -5mm \hspace{5mm}

Is there an analogue of Theorem~\ref{girth_l2} for embeddings into
$l_1$?
\begin{open}\rm
How small can $c_1(G)$ be for a a graph $G$ of girth $g$ in which
all vertices have degree $\ge 3$? Specifically, can $c_1(G)$ stay
bounded as $g$ tends to $\infty$?
\end{open}

\section{Ramsey-type theorems for metric spaces}

\vskip-5mm \hspace{5mm}

The philosophy of modern Ramsey Theory, (as developed e.g.
in~\cite{ramsey_book}) can be stated as follows: Large systems
necessarily contain substantial ``islands of order". Dvoretzky's
Theorem certainly falls into this circle of ideas. But what about
the metric analogues?
\begin{open}\rm
What is the largest $f(\cdot,\cdot)$ so that every $n$-point
metric $(X,d)$ has a subset $Y$ of cardinality $\ge f(n,t)$ with
$c_2(Y) \le t$? (We mean, of course, the metric $d$ restricted to
the set $Y$.)
\end{open}
For $t$ close to $1$, the answer is known, namely, $f(n,t) =
\Theta(\log n)$. For larger $t$ the behavior is known to be
different \cite{Ba:Li:Me:Na}.

\bibliographystyle{alpha}

\label{lastpage}

\end{document}